\documentclass[letterpaper, paper,10pt]{AAS}		% for final proceedings (20-page limit)
\usepackage[utf8]{inputenc}  % Helps a lot with arXiv submission!!!
\usepackage{bm}
\usepackage{amsmath}
\usepackage{subfigure}
\usepackage[colorlinks=true, pdfstartview=FitV, linkcolor=black, citecolor= black, urlcolor= black]{hyperref}
\usepackage{footnpag}			      	% make footnote symbols restart on each page

\usepackage{algorithm, algpseudocode}

\usepackage[superscript]{cite}
\usepackage{amssymb}
\usepackage{amsmath}
\usepackage{letltxmacro}
\LetLtxMacro{\originaleqref}{\eqref}
\renewcommand{\eqref}{Eq.~\originaleqref}
\newtheorem{remark}{Remark}

\newcommand{\Real}{\mathbb R}

\newcommand{\real}[1]{{\mathbb R}^{#1}}
\newcommand{\norm}[1]{\left\Vert#1\right\Vert}

\newcommand{\bb}{{\boldsymbol b}}

\newcommand{\be}{{\boldsymbol e}}
\newcommand{\bff}{{\boldsymbol f}}

\newcommand{\bh}{{\boldsymbol h}}

\newcommand{\bp}{{\boldsymbol p}}

\newcommand{\bu}{{\boldsymbol u}}

\newcommand{\bw}{{\boldsymbol w}}
\newcommand{\bx}{\boldsymbol x}

\newcommand{\bxf}{{\bx(\cdot)}}  % f for function
\newcommand{\buf}{{\bu(\cdot)}}  % f for function
\newcommand{\bB}{{\boldsymbol B}}

\newcommand{\bP}{{\boldsymbol P}}
\PaperNumber{25-689}
\begin{document}

\title{A MILLION-POINT FAST TRAJECTORY OPTIMIZATION SOLVER}
%----------------------------------------------------------------------------------
%
\author{A. Javeed$^*$ D. P. Kouri$^*$ D. Ridzal\thanks{Optimization and Uncertainty Quantification, Sandia National Laboratories, Albuquerque, NM, 87185.}\ \
J. D. Steinman\thanks{Graduate Research Assistant, Computation Applied Mathematics and Operations Research, Rice University, Houston, TX, 77005}
\ and I. M. Ross\thanks{Distinguished Professor and Program Director, Control and Optimization, MAE, Naval Postgraduate School, Monterey, CA, 93943.}
}
\maketitle{} 		
%==========================================================================
%
\begin{abstract}
One might argue that solving a trajectory optimization problem over a million grid points is preposterous.  How about solving such a problem at an incredibly fast computational time? On a small form-factor processor?  Algorithmic details that make possible this trifecta of breakthroughs are presented in this paper. The computational mathematics that deliver these advancements are: (i) a Birkhoff-theoretic discretization of optimal control problems, (ii) matrix-free linear algebra leveraging Krylov-subspace methods, and (iii) a near-perfect Birkhoff preconditioner that helps achieve $\mathcal{O}(1)$ iteration speed with respect to the grid size,~$N$. A key enabler of this high performance is the computation of Birkhoff matrix-vector products at $\mathcal{O}(N\log(N))$ time using fast Fourier transform techniques that eliminate traditional computational bottlenecks. A numerical demonstration of this unprecedented scale and speed is illustrated for a practical astrodynamics problem.

\end{abstract}

\section{Introduction}
\begin{figure}[h!]
      \centering
      {\parbox{0.9\columnwidth}{
      \centering
      {\includegraphics[width = 0.55\columnwidth]{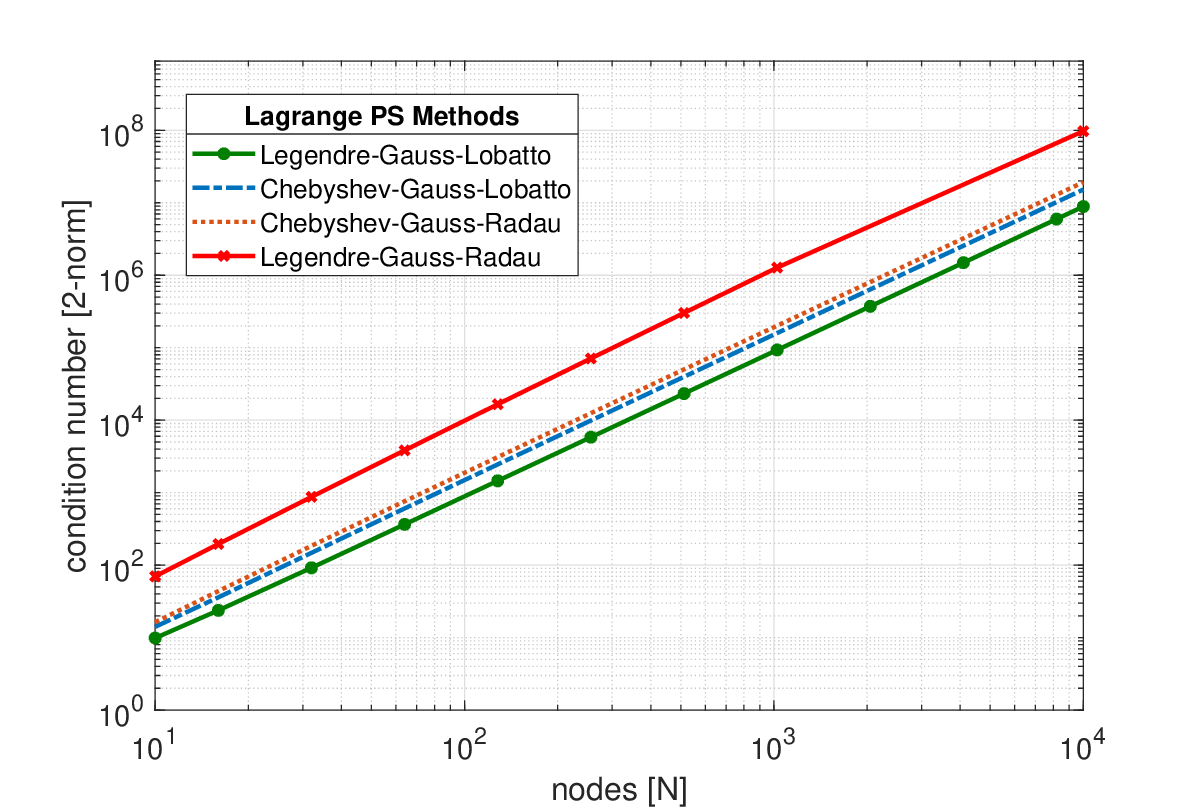}}
      \caption{Condition Numbers for $D$-Matrix-Based PS Methods\cite{fastmesh,newBirk-part-I}.}
      \label{fig:PScondNumbers}
      }
      }
\end{figure}

Consider the grand challenge of solving a nontrivial, nonlinear, nonconvex trajectory optimization problem over a million grid points. This challenge problem was posed in 2017 at a SIAM conference\cite{ross-million-2017} as an enabling technology to solve a new breed of orienteering\cite{ross-TSP-nolcos-2016} and traveling-salesperson problems\cite{ross-TSP-arXiv-2020} that are governed by differential equations. Nonetheless, it was argued in [\citen{ross-million-2017}] that if an ordinary trajectory optimization was solvable over a million grid points within reasonable space and time complexity bounds\cite{complexity-book-2007}, it would likely usher in a new era in computational optimal control. To put this point in perspective\cite{perspective}, we start at the beginning and consider the state of practice\cite{conway:survey,trelat:survey} in trajectory optimization methods as depicted in Figure~\ref{fig:TrajOptPerspective} (from Ref.~[\citenum{newBirk-part-I}]):
%
%======================================================================================
\begin{figure}[h!]
      \centering
      {\parbox{0.9\columnwidth}{
      \centering
      {\includegraphics[width = 0.65\columnwidth]{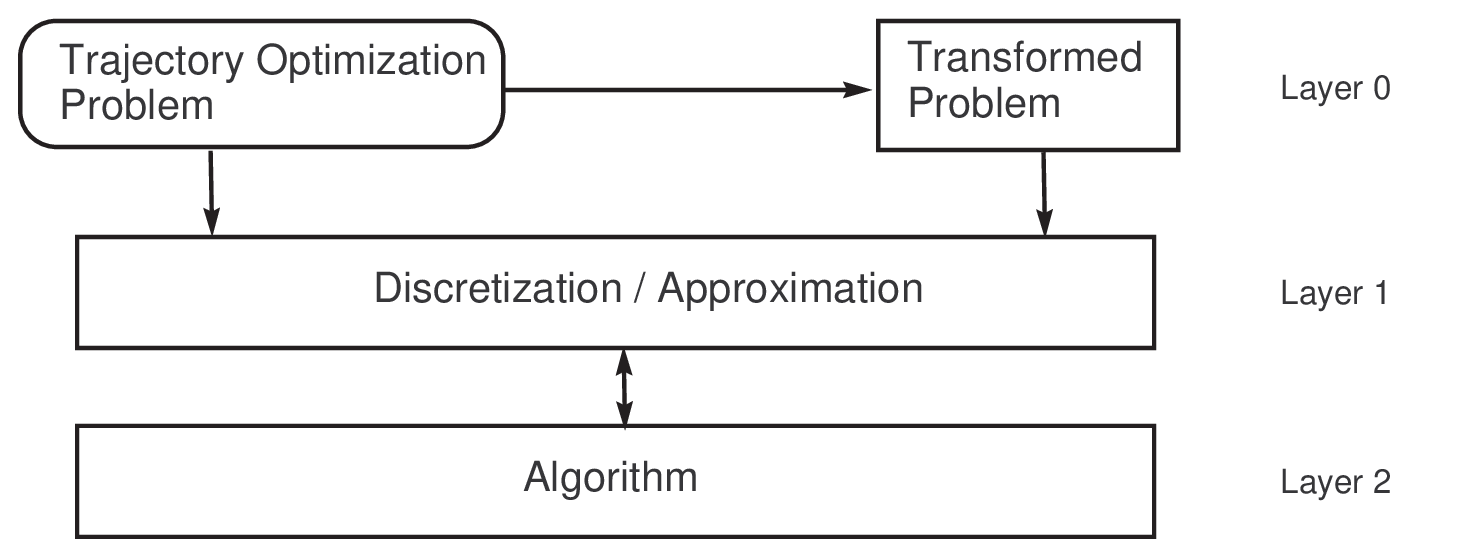}}
      \caption{A Schematic for the State-of-Practice in Trajectory Optimization Methods. From Reference~[\citenum{newBirk-part-I}].}
      \label{fig:TrajOptPerspective}
      }
      }
\end{figure}
%==========================================================================================
%
%
\begin{enumerate}
\item[($i$)] A decision is made to either work ``directly'' with the given problem or its transformation to its necessary conditions of optimality by an application of Pontryagin's Principle\cite{ross-book,longuski}. This is labeled Layer~0 in Figure~\ref{fig:TrajOptPerspective}. In the latter case, the process is known as an ``indirect'' method\cite{conway:survey,trelat:survey}.
\item[($ii$)] Regardless of the method being direct or indirect, a discretization is selected to convert the problem or its transformation into a continuous, finite-dimensional, optimization or a root-finding problem\cite{ross:roadmap-2005}. This is Layer~1 in Figure~\ref{fig:TrajOptPerspective}.
\item[($iii$)] The resulting finite-dimensional optimization or root-finding problem is passed to an off-the-shelf optimizer (like SNOPT\cite{snopt.paper} or IPOPT \cite{ipopt}) to be solved in Layer~$2$.
\end{enumerate}
Barring a few exceptions\cite{hr}, many optimizers make direct use of a Jacobian and/or a Hessian by way of sparsity patterns. That is, passing the optimizer only the matrix indices and values of the nonzero entries results in a more efficient linear algebra, including efficient linear system solves\cite{snopt.paper,ipopt}.  The bottleneck in these traditional ideas are space and time complexity\cite{complexity-book-2007,ross:Hessians}. Space complexity refers to the amount of computer memory required to solve a problem while time complexity refers to the total number of floating point operations required to achieve a solution within some given accuracy bounds\cite{ross:Hessians}. To better understand this bottleneck consider one of the most popular methods to solve trajectory optimization problems, namely a pseudospectral (PS) method\cite{PSReview-ARC-2012}.  The PS differentiation matrix\cite{boyd,trefethen-2000}, $D$, is a dense $N_n \times N_n$ matrix where $N_n := N+1$ is the number of node points or grid points.  If $N_n = 10^6$, then storing $D$ requires computer memory equivalent to storing $10^{12}$ numbers or $8$ terrabytes (TB).  Even if 8 TB were to be a small number at some point in the future, we argue that it is still more efficient and desirable to use, say, $8 \times 10^{-9}$ TB of memory or $0.008$ gigabytes (GB) or memory. The mechanism to achieve this feat is, in principle, quite simple: do not store $D$!  The differentiation matrix is an operator.  In trajectory optimization and other applications of PS methods\cite{boyd,trefethen-2000}, $D$ is not used by itself, rather what is important is the matrix-vector product $DX$, $X \in \real{N_n}$. Storing the product $DX$ over a million grid points takes only $0.008$ GB of memory. Thus, the ``trick'' is to compute $DX$ fast and without storing $D$.  This trick is well-known\cite{boyd,trefethen-2000}: compute $DX$ using a fast Fourier transform (FFT).  The FFT approach requires one to select a Chebyshev grid in order to exploit the connection between Chebyshev polynomials and the discrete cosine transform (DCT). Or simply put, one must use the Chebyshev PS method for trajectory optimization\cite{fahroo:cheb-jgcd,cheb-costate}.  In principle, using an FFT also solves the time complexity problem.  If the product $DX$ is computed using standard matrix multiplication, the time complexity of this process is $N_n^2$ multiplications which amounts to $10^{12}$ floating point operations over a million grid points.  Because an FFT requires only $\mathcal{O}(N\log(N))$ multiplications, the product $DX$ can be computed in approximately $14\times 10^6$ floating point operations over a million grid points, or about $72,000$ times faster than standard matrix multiplication!  Although the preceding ideas have been known for quite some time\cite{boyd,trefethen-2000} it has never been used to solve a million-point trajectory optimization problem using a Chebyshev PS method because of the well-known problem of large condition numbers\cite{boyd,trefethen-2000,fastmesh,newBirk-part-I,newBirk-part-II} associated with differentiation matrices (see Figure~\ref{fig:PScondNumbers}).
%
%======================================================================================
%==========================================================================================
%
That is, the condition number of a $D$-matrix-based PS method (i.e., a Lagrange interpolation-based method) increases at the rate of $\mathcal{O}(N^2)$ no matter the selection of the Gaussian grid.  As shown in Figure~\ref{fig:PScondNumbers}, the best and worst performers are the Legendre-Gauss-Lobatto (LGL) grid and the Legendre-Gauss-Radau grid respectively.  The Chebyshev-Gauss-Lobatto (CGL) grid is close to the best performer, namely, LGL, but eventually all grid points suffer from excessively large condition numbers. A widely-used approach for overcoming this problem in trajectory optimization is to use a concatenation of low-order polynomials, namely a PS knotting method\cite{knots,auto-knots,spec-alg} which is the time adaptation of an $hp$-adaptive\cite{hp-adaptive} or $hp$-element method\cite{spectral-hp,boyd} used in the context of spatial discretizations.  Although this simple idea solves the condition number problem, it comes at a hefty price of exponentially reduced convergence rate\cite{kang-rate-2010,Kang_2008_convergence,PSReview-ARC-2012} in the case of trajectory optimization problems.  Additional computation cost is incurred in determining the number and location of the knot points\cite{knots,auto-knots} (i.e., ``$h$-segments'').  For instance, in Reference~[\citen{spec-alg}] the number of segments are determined by differentiating the control function to determine the locations of sharp changes (e.g. switches) in the trajectory.  PS knots are placed near these locations and the problem is resolved. To enhance efficiency, the number of node points (i.e., ``$p$-order'') in between the knot points is re-estimated\cite{auto-knots}.   This re-estimation is performed using Jackson's theorem\cite{spec-alg,auto-knots} and other convergence criteria\cite{kang-rate-2010,Kang_2008_convergence,PSReview-ARC-2012}. After the problem is re-solved, the entire process is repeated to add/remove/relocate all the knot points\cite{knots,auto-knots}.

The introduction of the universal Birkhoff theory for trajectory optimization\cite{newBirk-part-I,newBirk-part-II} eliminated all of these problems in one fell swoop.
%
%======================================================================================
\begin{figure}[h!]
      \centering
      {\parbox{0.9\columnwidth}{
      \centering
      {\includegraphics[width = 0.55\columnwidth]{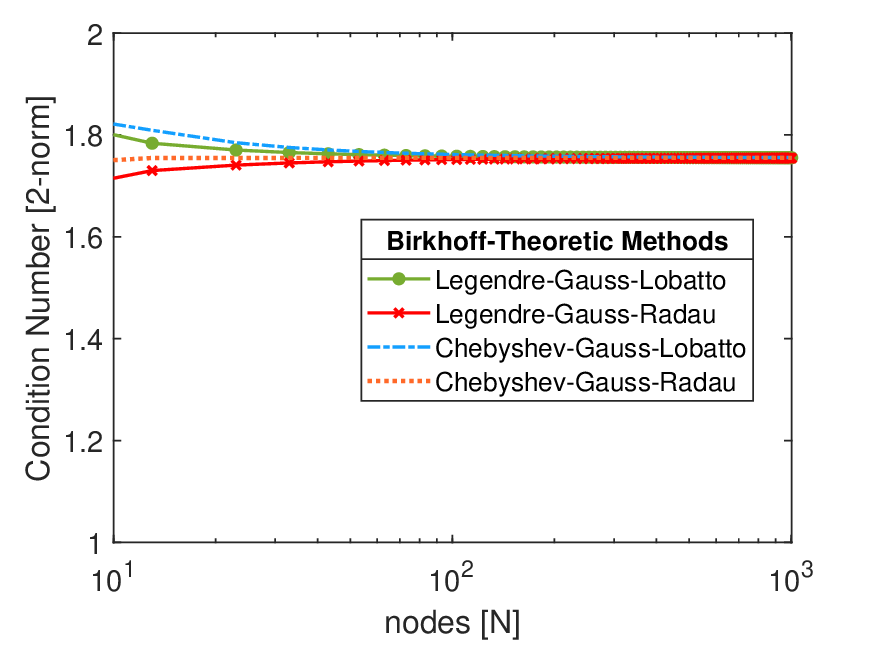}}
      \caption{Condition Numbers for Birkhoff-Theoretic Methods\cite{newBirk-part-I}.}
      \label{fig:BirkCondNumbers}
      }
      }
\end{figure}
%==========================================================================================
%
In a Birkhoff-theoretic method, the conditions numbers are small and astonishingly ``flat;'' i.e., vary as $\mathcal{O}(1)$.  See Figure~\ref{fig:BirkCondNumbers}.
Furthermore, a Birkhoff-theoretic method retains all the desirable properties of a standard PS method such as exponential convergence and computation over a Chebyshev grid\cite{newBirk-part-I,newBirk-part-II}. In other words, a Birkhoff-Chebyshev method sets the stage for solving a million-point trajectory optimization problem. It will be apparent shortly that a Birkhoff-Chebyshev method is simply a starting point.  A substantial set of additional ideas and algorithms are necessary to meet the challenges of a fast and efficient million-point trajectory optimization solver.  A quick summary of these additional ideas are:
\begin{enumerate}
\item The use of a DCT to generate a Birkhoff-Chebyshev matrix-vector product.
\item The construction of a lower-triangular Birkhoff-like matrix that enables its use as part of an efficient preconditioner for a Newton/Karush-Kuhn-Tucker (KKT) system\cite{ross:Hessians}.
\item The use of preconditioned matrix-free linear algebra\cite{matrix-free-ref-1,matrix-free-ref-2,hr} to solve the KKT system.
\end{enumerate}
From these requirements, it is clear why current off-the-shelf optimization packages like SNOPT\cite{snopt.paper} and IPOPT\cite{ipopt} are unsuitable as an engine for matrix-free trajectory optimization. Stated differently, it is necessary to build from the ground up an entirely new integrated approach to produce a fast million-point trajectory optimization solver .

\section{An Overview of the Million-Point Trajectory Optimization Solver}
It is not possible to simply use FFTs and a Chebyshev-Birkhoff discretization to solve a million point trajectory optimization problem using off-the-shelf optimization routines. This is because (as noted in the Introduction section) many such routines\cite{snopt.paper,ipopt} require explicit Jacobian and/or Hessian information.  Storing a Jacobian and/or Hessian matrix runs up against the bottlenecks of space and time complexity issues noted previously\cite{ross:Hessians}. These bottlenecks can be overcome by using the concept of matrix-free methods\cite{matrix-free-ref-1,matrix-free-ref-2,hr}.  That is, instead of using explicit Jacobian/Hessian information, we use the product of these matrices with a vector without computing the matrix.  This is the key idea behind matrix-free linear algebra\cite{matrix-free-ref-1,matrix-free-ref-2,hr}.

\subsection{A Brief Introduction to Matrix-Free Linear Algebra}
Consider the fundamental problem of solving the linear matrix equation,
\begin{equation}\label{eq:Ax=b}
A\,\chi = \beta
\end{equation}
where $A$ is an $n \times n$ dense square matrix of full rank, $\chi \in \real{n}$ and $\beta \in \real{n}$.  A standard ``direct method''\cite{trefethen-bau-1997,golub-book-2013} for solving \eqref{eq:Ax=b} is to factorize $A$ (by an $LU$ decomposition, for example) and sequentially solve the factored system.  The space- and time complexity of such a method is $\mathcal{O}(n^2)$ and approximately $\mathcal{O}(n^3)$ respectively\cite{trefethen-bau-1997,golub-book-2013}. The space complexity of $O(n^2)$ was already discussed in the Introduction section.  The time complexity of $O(n^3)$ is approximate because the best method of solving \eqref{eq:Ax=b} is currently\cite{matrix-free-ref-1} at $O(n^{\omega})$, where $\omega \approx 2.3$. Taking $O(n^3)$ as a baseline case for solving \eqref{eq:Ax=b} for $n=10^6$ takes approximately $11$ days of compute time on a terraflop processor\cite{ross:Hessians}. Recall that storing $A$ requires $8$ TB of memory. Stated differently, even if FFTs and Birkhoff-Chebyshev discretizations were used in Layer~1 of the trajectory optimization process shown in Figure~\ref{fig:TrajOptPerspective}, all of its advantages would be completely lost if direct methods were used in Layer~2.  In other words, using popular solvers like SNOPT and IPOPT would defeat the entire process because of their explicit use of Jacobian and/or Hessian information. In the world of matrix-free linear algebra, \eqref{eq:Ax=b} is solved ``indirectly\footnote{Note that the words ``direct'' and ``indirect'' as commonly used in linear algebra are not connected in any way to similar verbiage used in trajectory optimization\cite{perspective,conway:survey,trelat:survey}.}'' via an iterative process\cite{saad.2003,CTKelley.1999}.  This approach does not require $A$; instead, it relies on the (efficient) computation of the product $A\chi$ to solve \eqref{eq:Ax=b}.  The basic idea behind a matrix-free method is quite simple as described by Algorithm~1:\\[1em]
%
%{\centering{%
\fbox{%
    \parbox{0.95\linewidth}{%
\textbf{Algorithm 1: Rudimentary matrix-free method for solving \eqref{eq:Ax=b}}:
\begin{enumerate}
\item Let $\chi^0$ be a guess to the solution of \eqref{eq:Ax=b}.
\item Compute the product $A\chi^0$ (without using $A$ explicitly such as through the use of FFTs).
\item If $\norm{A\chi^0-\beta} \le \epsilon$, where $\epsilon > 0$ is a given solution tolerance, stop and exit.  Otherwise continue.
\item For $k = 0, 1, \ldots$ compute $\chi^{k+1} = \chi^k + \Delta \chi^k$ where $\Delta \chi^k$ is found iteratively using matrix-vector products of $A\chi^k$ (without using $A$ explicitly).
\item Stop the iteration when $\norm{A\chi^k-\beta} \le \epsilon$
\end{enumerate}
}%
}
%}}
\\[1em]
%------------
\begin{remark}
If the product $A\chi$ in the preceding algorithm is computed by using $A$ explicitly then the resulting process is an iterative method for solving \eqref{eq:Ax=b} rather than a matrix-free method. See [\citen{gmres.MATLAB}] for an illustrative example on how to convert the product of a sparse matrix $A$ with a vector $\chi$ to produce a matrix-free computation of $A\chi$.
\end{remark}
%------
Deferring a discussion on the key step of computing $\Delta \chi^k$, it is obvious that the complexity of solving \eqref{eq:Ax=b} using a matrix-free approach depends upon the efficiency of computing $A\chi$ and the number of iterations it takes to generate $\norm{A\chi^k-\beta} \le \epsilon$.  Both the accuracy and efficiency (i.e., number of iterations) of this approach depends strongly\cite{CTKelley.1999,trefethen-bau-1997,golub-book-2013} on the condition number of $A$. In trajectory optimization, this $A$-matrix comprises\cite{ross:Hessians} discretized Jacobians and Hessians and hence it is necessary to never store this ``inner-loop'' matrix so that \eqref{eq:Ax=b} can be solved quickly and accurately over a million grid points.

\subsection{Introducing a Matrix-Free Trajectory Optimization Solver}
From the discussions of the preceding subsection, it is clear that constructing a million point trajectory optimization solver requires matrix-free linear algebra. A popular matrix-free method for solving \eqref{eq:Ax=b} is a conjugate gradient (CG) method\cite{saad.2003,CTKelley.1999}. Interestingly, many solvers (e.g. SNOPT) offer a conjugate gradient option to solve \eqref{eq:Ax=b}.  As noted in the previous subsection, the most efficient way to use a matrix-free method, and hence a CG method, is to never compute or store~$A$. Unfortunately, solvers like SNOPT implement a CG algorithm by asking a user to explicitly provide the Jacobian information (or its sparsity pattern for internal computation).  In a Birkhoff-theoretic method, the $B$-matrix is dense\cite{newBirk-2023}; hence, the advantage offered by a CG method is quickly lost if matrix-vector products are computed using ordinary matrix multiplication methods. Thus, the minimum requirements for a million point trajectory optimization solver are:
\begin{enumerate}
\item Matrix-free computation of Jacobian and Hessian matrix-vector products; and
\item Efficient incorporation of problem-specific preconditioners.
\end{enumerate}
To the best of the authors' knowledge, such trajectory optimization routines are currently unavailable.  As a result, we construct a first-of-a-kind method that incorporates matrix-free linear algebra with problem-specific preconditioners.  The details of this new approach are discussed in the sections to follow.

\section{Birkhoff Discretizations of Generic Optimal Control Problems}
As noted in the Introduction section, the challenge posed in [\citen{ross-million-2017}] was to solve a generic, nonlinear, nonconvex optimal control problem over a million grid points.  We pose this generic problem as\cite{ross-book}:
%
%=====================================================================================
\begin{eqnarray}
&\bx(t) \in \real{N_x}, \quad \bu(t) \in \real{N_u}, \quad t \in [t^a, t^b], \quad \bp \in \real{N_p}   & \nonumber\\
& (\textsf{$\bP$}) \left\{
\begin{array}{lrl}
\emph{Minimize } & J[\bx(\cdot), \bu(\cdot), t^a, t^b, \bp] =& E(\bx(t^a),\bx(t^b), t^a, t^b, \bp)\\
& & \quad \displaystyle + \int_{t^a}^{t^b} F(\bx(t), \bu(t), t, \bp)\, dt \\
\emph{Subject to}& \dot\bx(t) =& \bff(\bx(t), \bu(t), t, \bp)  \\
& \be^L \le & \be(\bx(t^a), \bx(t^b), t^a, t^b, \bp) \le \be^U \\
& \bh^L \le & \bh(\bx(t), \bu(t), t, \bp) \le \bh^U
\end{array} \right. & \label{eq:probPbold}
\end{eqnarray}
%=====================================================================================
%
where, the symbols imply the following:
\begin{itemize}
\item $N_{(\cdot)} \in \mathbb{N}$ represents the number of variables or functions implied by the subscript $(\cdot)$;
\item $\bx \in \real{N_x}$ is the state variable, $\bu \in \real{N_u}$ is the control variable, $t \in \Real$ is the clock time, $t^a \in \Real$ is an initial clock time, $t^b > t^a$ is a final clock time, and $\bp \in \real{N_p}$ is an optimization parameter;
\item $(\bxf, \buf)$ is the state-control function pair (i.e., an element in some Sobolev space\cite{ross-book,vinter,clarke-2013book});
\item $J: \big(\bx(\cdot), \bu(\cdot), t_0, t_f, \bp\big) \mapsto \Real $ is the standard (``Bolza'') cost functional, $E: \real{N_x} \times \Real \times \real{N_x} \times \Real \times \real{N_p} \to \Real $ is the endpoint (``Mayer'') cost and $F: \real{N_x} \times \real{N_u} \times \Real \times \real{N_p} \to \Real$ is the running (``Lagrange'') cost\cite{longuski,ross-book};
\item $\bff$ is the dynamics function, $\bff:\real{N_x} \times \real{N_u} \times \Real \times \real{N_p} \to \real{N_x} $;
\item $\be$ is the endpoint constraint function, $\be: \real{N_x} \times \real{N_x} \times \Real \times \Real \times \real{N_p} \to \real{N_e} $ whose values are bounded by $[\be^L, \be^U]$, $\be^L \le \be^U$; and,
\item $\bh$ is the mixed state-control path constraint function, $\bh:\real{N_x} \times \real{N_u} \times \Real \times \real{N_p} \to \real{N_h} $ whose values are bounded by $[\bh^L, \bh^U]$, $\bh^L \le \bh^U$.
\end{itemize}
Problem~$(\bP)$ constitutes the left box of Layer~0 in Figure~\ref{fig:TrajOptPerspective}.  The production of Layer~1 for this problem is described in detail in References~[\citenum{newBirk-part-I}] and [\citen{newBirk-part-II}].  In the interest of completeness, we review this process using the fundamentals of the universal Birkhoff theory.\cite{newBirk-part-I,newBirk-2023}

\subsection{A Brief Introduction to the Universal Birkhoff Interpolants}
\begin{figure}[h!]
      \centering
      {\parbox{0.9\columnwidth}{
      \centering
      {\includegraphics[width = 0.65\columnwidth]{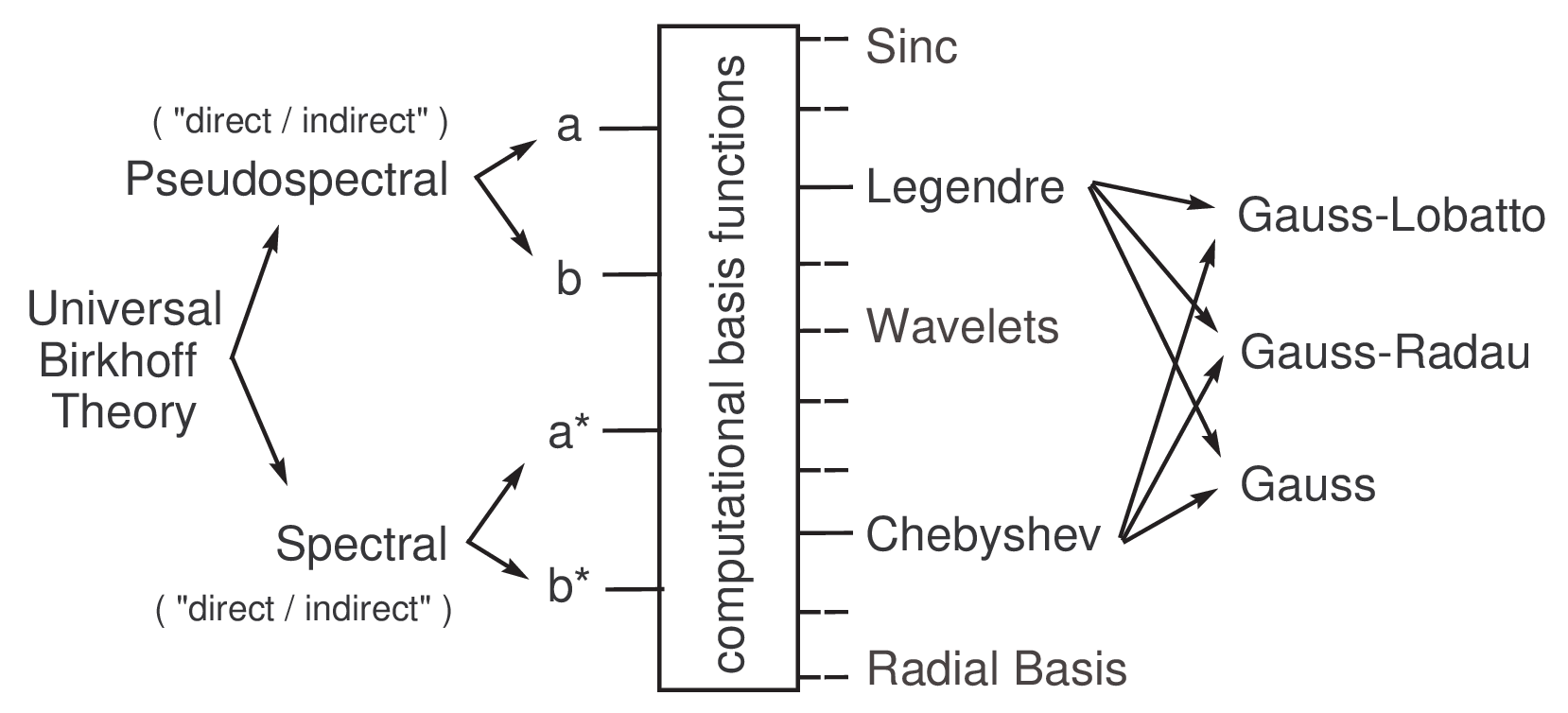}}
      \caption{Variants of Birkhoff Discretization Methods for (Direct/Indirect) Trajectory Optimization Methods; see Layer~1 in Figure~\ref{fig:TrajOptPerspective}. This figure is adapted from Reference~[\citenum{DIDO:arXiv}].}\label{fig:BirkhoffPStypes}
      }
      }
\end{figure}

As introduced in References~[\citen{newBirk-part-I}] and [\citen{newBirk-2023}], there are a large number of variants of the universal Birkhoff interpolants leading to a wide variety of discretization methods and selection of basis functions and grid points.  The richness of this approach is illustrated in Figure~\ref{fig:BirkhoffPStypes}. Because we are largely concerned with the Chebyshev variant, we hereafter limit the description of the Birkhoff theory to the use of Chebyshev polynomials (of the first kind\cite{boyd,trefethen-2000}). In particular, we use the CGL route that is evident in Figure~\ref{fig:BirkhoffPStypes}.
%
%======================================================================================
%==========================================================================================
%

A CGL grid is defined on the interval $[-1, 1]$ and given in closed form by\cite{boyd,trefethen-2000},
\begin{equation}\label{eq:CGLpts}
\tau_j = -\cos\left(\frac{j \pi}{N} \right) = \cos\left(\frac{(N-j) \pi}{N} \right) , \quad j = 0, \ldots, N
\end{equation}
As has been noted several times before\cite{newBirk-part-I,newBirk-part-II,newBirk-2023,cheb-costate,fahroo:cheb-jgcd}, in sharp contrast to Legendre grids (be it Gauss, Radau or Lobatto), all Chebyshev grids can be generated using a single line of code.  Let $[-1, 1] \ni \tau \mapsto x(\tau) \in \Real$ be a continuously differentiable function.  To maintain consistency with the notation presented in [\citen{newBirk-part-I}] we let $\tau^a = -1 = -\tau^b$. A universal Birkhoff method over a CGL grid allows $\tau \mapsto x(\tau)$ to be approximated by two equivalent $a$- and $b$-expansions:
\begin{align}\label{eq:xN=Birka+b}
x^N(\tau) &:= x^a B_0^0(\tau) + \sum_{j=0}^N v_j B_j^a(\tau)
& x^N(\tau) := \sum_{j=0}^N v_j B_j^b(\tau) + x^b B_N^N(\tau)
\end{align}
where $(x^a, x^b) = (x(\tau^a), x(\tau^b))$, $B_0^0(\tau)$, $B^N_N(\tau)$, $ B_j^a(\tau)$ and $B_j^b(\tau)$, $j = 0, \ldots N$ are the Birkhoff basis functions that satisfy the following interpolating conditions\cite{newBirk-2023,newBirk-part-I}:
\begin{align}\label{eq:birk-conditions}
\begin{aligned}
B_0^0(\tau^a)           &= 1,            &&    B_N^N(\tau^b)   = 1,                      \\
 {d_\tau B}_0^0(\tau_i)   &=  0,          &&   {d_\tau B}_N^N(\tau_i)   = 0,   & i = 0, \ldots, N \\
B_j^a(\tau^a)           &=  0,            && B_j^b(\tau^b) = 0,                       & j = 0, \ldots, N \\
{d_\tau B}_j^a(\tau_i) &=  \delta_{ij},  && {d_\tau B}_j^b(\tau_i) =  \delta_{ij}, &   i = 0, \ldots, N,  & j = 0, \ldots, N
\end{aligned}
\end{align}
In \eqref{eq:birk-conditions}, $d_\tau$ is the derivative operator, $d/d\tau$, and  $\delta_{ij}$ is the Kronecker delta.  The variables $v_j, j = 0, \ldots, N$ in \eqref{eq:xN=Birka+b} are called the virtual variables or the virtual control variables\cite{DIDO:arXiv,newBirk-part-I}.

\subsection{A Simplified Presentation of Birkhoff Discretizations of Optimal Control Problems}
Using \eqref{eq:xN=Birka+b} to approximate Problem~$(\bP)$ generates its Birkhoff discretization/approximation.  This process generates Layer~1 in Figure~\ref{fig:TrajOptPerspective}. A complete description of this process and the final result are presented in [\citen{newBirk-part-II}]. For the purposes of completeness and illustrating this process we consider a ``distilled'' version\cite{newBirk-part-I} of Problem~$(\bP)$ given by,
%
%===============================================
\begin{eqnarray}
&x(t) \in \Real, \quad u(t) \in \Real, \quad \tau \in [\tau^a, \tau^b] = [-1, 1]   & \nonumber\\
& (P) \left\{
\begin{array}{lrl}
\emph{Minimize } && J[x(\cdot), u(\cdot)] := E(x(\tau^a), x(\tau^b))\\[0.5em]
\emph{Subject to}&& \dot x(\tau) = f(x(\tau), u(\tau)) \\
                 && e(x(\tau^a), x(\tau^b)) =  0
\end{array} \right. & \label{eq:probP}
\end{eqnarray}
%===============================================
%
%
%
Let, $x^N(\tau_i) = x_i$ and $u(\tau_i) = u_i$, $i = 0, \ldots, N$.  Define\cite{newBirk-2023,newBirk-part-I},
\begin{equation}\label{eq:manyDefs}
\begin{aligned}
X & :=  [x_0, \ldots, x_N]^T
& V &:=  [v_0, \ldots, v_N]^T
& U &:=  [u_0, \ldots, u_N]^T \\
f(X, U) &:=  [f(x_0, u_0), \ldots, f(x_N, u_N)]^T
& \bw_B &=  [w^{CC}_0, \ldots, w^{CC}_N]^T
& \bb &:=  [1, \ldots, 1]^T \\
[\bB^a]_{ij} & :=  [B^a_j({\tau_i})]
&[\bB^b]_{ij} &:=  [B^b_j({\tau_i})]
& N_n &:=  N+1
\end{aligned}
\end{equation}
where $w^{CC}_i, i = 0, \ldots, N$ are the Clenshaw-Curtis quadrature weights\cite{newBirk-2023} (for CGL points).  Then, the $a$-version of the Birkhoff discretization of Problem~$(P)$ is given by\cite{newBirk-part-I},
%
%===========================================================
\begin{align}
X \in \real{N_n}, \quad U \in \real{N_n}, \quad V \in \real{N_n}, \quad (x^a, x^b) \in \real{2}   & \nonumber\\
 \left(\textsf{$P_a^N$}\right) \left\{
\begin{array}{lrl}
\emph{Minimize } & J^N_a[X, U, V, x^a, x^b] :=& E(x^a,x^b)\\
\emph{Subject to} & X  = &x^a\,\bb + \bB^a V \\
& V =& f(X, U)  \\
& x^b = & x^a + \bw_B^T V \\
& e(x^a, x^b)  = & 0
\end{array} \right. & \label{eq:Prob-PNa}
\end{align}
%===============================================
%

\section{Fast Computation of Birkhoff Matrix-Vector Products }
Before discussing a method to solve Problem~$(P^N_a)$, we consider the more fundamental problem of computing the cost and constraint functions. We follow [\citen{newBirk-part-I,newBirk-part-II}] and assume, for the moment, that the functions $E, e$ and $f$ are zero-order oracles; i.e., they return only the values of the functions and not their gradients or Hessians. Then the computational cost associated with computing $E$ or $e$ is $\mathcal{O}(1)$, where the $\mathcal{O}(1)$ valuation is in terms of function calls and not the total computational cost (which would be problem dependent).  The computational cost associated with computing $f(X, U)$ is $\mathcal{O}(N)$.  These principles and ideas carry over to the general case of Problem~$(\bP)$; see [\citen{newBirk-part-II}] for details. From these conclusions and examining \eqref{eq:Prob-PNa} it is apparent that the computational cost associated with computing the cost and constraint functions of Problem~$(P^N_a)$ is dominated by the matrix-vector product $\bB^a V$.  As noted in the Introduction section, computing $\bB^a$ first and then performing the matrix-vector product is expensive and costs $\mathcal{O}(N^2)$ multiplications (and $\mathcal{O}(N^2)$ storage).  We now show that this expense can be dramatically reduced to $\mathcal{O}(N \log(N))$ through the use of an FFT (implemented by a DCT).  These details are described in the following subsections.

\subsection{Fast Matrix-Free Nodal to Modal to Nodal Transformations}
Differentiating the $a$-version of \eqref{eq:xN=Birka+b} (with respect to $\tau$) we get,
\begin{equation}\label{eq:xdot=v=sumT}
\dot x^N(\tau) := v(\tau) = \sum_{j=0}^N v_j \dot B_j^a(\tau) = \sum_{k=0}^N a_k T_k(\tau)
\end{equation}
where the last equality is just an (orthogonal) expansion of $v(\tau)$ (called a virtual variable in [\citen{DIDO:arXiv}]) in terms of Chebyshev polynomials $T_k(\tau), k = 0, \ldots, N$.  From \eqref{eq:birk-conditions} and \eqref{eq:xdot=v=sumT} we get a connection between the nodal values of $v(\tau)$ at $\tau_i, i = 0, \ldots, N$ and its Chebyshev modal coefficients $a_k, k = 0, \ldots, N$ given by,
\begin{equation}\label{eq:vi=sumT}
v(\tau_i) =  v_i  = \sum_{k=0}^N a_k T_k(\tau_i), \quad i = 0, \ldots, N
\end{equation}
Equation~(\ref{eq:vi=sumT}) can be written as a modal-to-nodal (or coefficient-to-value) linear transformation according to,
\begin{equation}\label{eq:V=Ta}
V = \left[
              \begin{array}{c}
                v_0 \\
                \vdots \\
                v_N \\
              \end{array}
            \right] = \underbrace{\left[
      \begin{array}{ccc}
        T_0(\tau_0) & \cdots & T_N(\tau_0) \\
        \vdots & \ddots & \vdots \\
        T_0(\tau_N) & \cdots & T_N(\tau_N) \\
      \end{array}
    \right]}_{\mathbf{T}} \left[
              \begin{array}{c}
                a_0 \\
                \vdots \\
                a_N \\
              \end{array}
            \right]
\end{equation}
As shown in Appendix~A (see \eqref{eq:T=cos}) the Chebyshev polynomials can be written as,
\begin{equation}
T_k(\tau) = \cos(k\theta(\tau)), \quad \theta(\tau) = \cos^{-1}(\tau), \quad \tau \in [-1, 1], \quad k = 0, \ldots, N
\end{equation}
Thus, the value of $T_k(\tau)$ at the CGL point $\tau_j$ (see \eqref{eq:CGLpts}) is given by,
\begin{equation}\label{eq:Tkj=coskj}
T_k(\tau_j) = \cos\left(\frac{k(N-j) \pi}{N} \right), \quad j = 0, \ldots, N
\end{equation}
Hence, each element of the $(N+1)\times(N+1)$ matrix $\mathbf{T}$ in \eqref{eq:V=Ta} can be computed as fast as the computation of a cosine function.  It thus seems obvious that $V$ can be easily computed by determining $\mathbf{T}$ first and then computing the matrix-vector product indicated in \eqref{eq:V=Ta}.  However, as noted previously\cite{ross:Hessians}, this process is quite expensive in terms of space and time complexity (i.e., $\mathcal{O}(N^2)$).  \emph{The key to performing a more efficient computation of the transformation indicated in \eqref{eq:V=Ta} (and its inverse) is to never compute/store the matrix $\mathbf{T}$.  This is the essential idea behind matrix-free methods.}

The matrix-free version of \eqref{eq:V=Ta} relies on rewriting \eqref{eq:vi=sumT} in the form of a DCT and using an FFT-type efficiency to perform the requisite operations. Although this matrix-free approach is well-known\cite{boyd,trefethen-2000}, we briefly describe it here in order to contrast it with the direct matrix-vector product implied by \eqref{eq:V=Ta}.

Substituting \eqref{eq:Tkj=coskj} in \eqref{eq:vi=sumT} we get,
\begin{multline}\label{eq:vm=almostDCT-1}
v_i  = \sum_{k=0}^N a_k T_k(\tau_i) =  \sum_{k=0}^N a_k \cos\left(\frac{k(N-i) \pi}{N} \right)\\
= a_0 +  \sum_{k=1}^{N-1} a_k \cos\left(\frac{k(N-i) \pi}{N} \right) + (-1)^{(N-i)} a_N
\end{multline}
Equation~(\ref{eq:vm=almostDCT-1}) is almost the same as the type-I DCT of real numbers $\widehat{a}_0, \ldots, \widehat{a}_N$, given by,
\begin{equation}\label{eq:vm=DCT-1}
v_j^{DCT}  = \frac{\widehat{a}_0}{2} +  \sum_{k=1}^{N-1} \widehat{a}_k \cos\left(\frac{k\, j\, \pi}{N} \right) + (-1)^{j}\ \frac{\widehat{a}_N}{2}, \quad j = 0, \ldots, N
\end{equation}
Hence, we rewrite \eqref{eq:vm=almostDCT-1} by setting $N-i =j$ to produce,
\begin{equation}\label{eq:coeff2v-1}
v_{(N-j)} = \frac{2\,a_0}{2} +  \sum_{k=1}^{N-1} a_k \cos\left(\frac{k\,j \pi}{N} \right) + (-1)^{j} \frac{2\, a_N}{2}
\end{equation}
Equation~(\ref{eq:coeff2v-1}) generates a  matrix-free version of \eqref{eq:V=Ta} as follows:\\[1em]
%
%=========================================================
\fbox{%
    \parbox{0.95\linewidth}{%
\textbf{\texttt{Fasta2v}: A fast algorithm for computing the Chebyshev nodal values from its modal coefficients}:
\begin{enumerate}
\item Given $a_0, a_1, \ldots, a_N$, multiply the first and last term by two.
\item Apply DCT-I to $2 a_0, a_1, \ldots, 2 a_N$ to generate $v_N, v_{N-1}, \ldots, v_0$.
\item Rearrange the result according to $v_0, v_1, \ldots, v_N$.
\end{enumerate}
}}\\ \\
%===============================================
%
Obviously, Steps~1 and 3 are computationally inexpensive.  Step~2 can be performed in $\mathcal{O}(N\log(N))$ operations.  %as follows:
Note also that the space complexity of using a DCT is just $\mathcal{O}(N)$. This dramatic reduction in space and time complexity also holds for the inverse problem of generating $a_0, a_1, \ldots, a_N$ from $v_0, v_1, \ldots, v_N$.  The naive approach is, of course, to simply invert $\mathbf{T}$ which is extremely expensive (i.e., $\mathcal{O}(N^3)$).  The fast matrix-free approach is to reverse the three steps used in computing \eqref{eq:V=Ta}.  Although this point is quite obvious, it is derived from first principles in Appendix~B.

We now use these principles to show how to compute $\bB^a V$ without computing $\bB^a$ first.

\subsection{An FFT Algorithm for Computing $\bB^a V$}
\begin{figure}[h!]
      \centering
      {\parbox{0.9\columnwidth}{
      \centering
      {\includegraphics[width = 0.5\columnwidth]{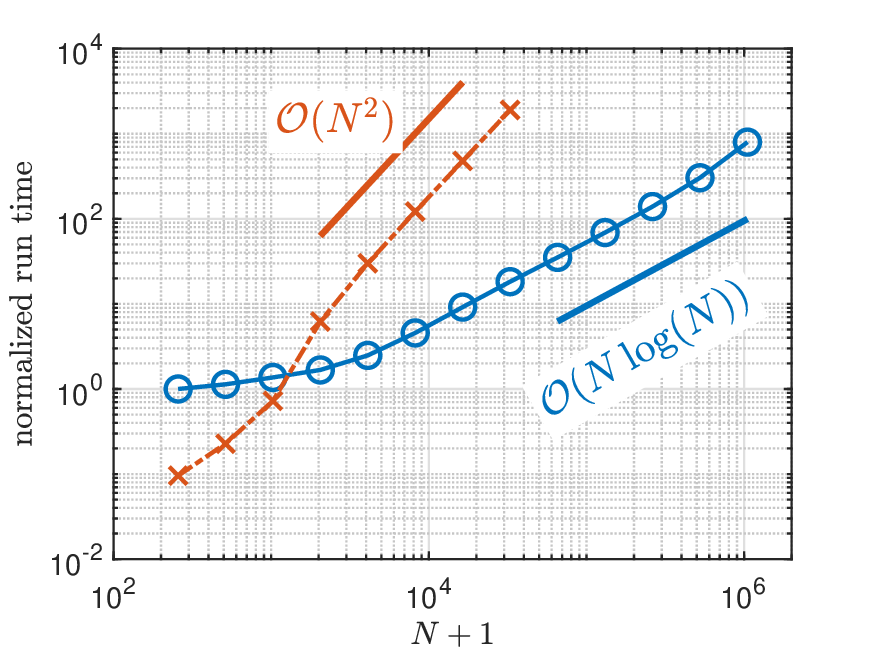}}
      \caption{Relative Run Times for Computing $\bB^aV$ Using \textbf{\texttt{FastBV}} (blue) and as Ordinary Matrix Multiplication (red).}
      \label{fig:BVisNlogN}
      }
      }
\end{figure}

The Birkhoff matrices $\bB^a$ and $\bB^b$ comprise values of the Birkhoff interpolants at the node points (see \eqref{eq:xN=Birka+b} and \eqref{eq:manyDefs}).  Hence, the $i^{th}$-row of $\bB^a V$ (denoted as $[\bB^a V]_i$) may be viewed equivalently as
\begin{equation}\label{eq:BVi}
\left[\bB^a V\right]_i := \sum_{j=0}^N v_j B_j^a(\tau_i) = \sum_{j=0}^N v_j B_j^a(\tau){\big\vert}_{\tau = \tau_i}
\end{equation}
The nuanced idea behind the last equality in \eqref{eq:BVi} is that it commutes the evaluation operation $\tau=\tau_i$ with the computation of the Birkhoff interpolant.  That is, the second equality suggests that that $\sum_{j=0}^N v_j B_j^a(\tau)$ may be computed first (for any $\tau$) and then evaluated at $\tau=\tau_i$ afterwards.  Taking this idea further, we write\cite{newBirk-2023},
\begin{equation}\label{eq:BVt-1}
\sum_{j=0}^N v_j B_j^a(\tau) = \sum_{j=0}^N v_j \int_{-1}^{\tau}\dot B_j^a(\xi)d\xi = \sum_{k=0}^N a_k \int_{-1}^{\tau} T_k(\xi) d\xi
\end{equation}
where the last equality follows from \eqref{eq:xdot=v=sumT}. It is shown in [\citen{newBirk-2023}] and the references contained therein, that
\begin{subequations}\label{eq:intChebPoly=}
\begin{align}
\int_{-1}^\tau T_k(\xi)d\xi &= \frac{T_{k+1}(\tau)}{2(k+1)} - \frac{T_{k-1}(\tau)}{2(k-1)} - \frac{(-1)^k T_0(\tau)}{k^2-1}, \quad \text{for } \quad k \ge 2\\
\int_{-1}^\tau T_1(\xi)d\xi &= \frac{T_2(\tau)} {2(2)}- \frac{T_0 (\tau)}{2(2)} \\
\int_{-1}^\tau T_0(\xi)d\xi &= T_1(\tau) + T_0(\tau)
\end{align}
\end{subequations}
Hence, \eqref{eq:BVt-1} may be written as,
\begin{multline}\label{eq:BVt-2}
\sum_{j=0}^N v_j B_j^a(\tau)  = \sum_{k=0}^N a_k \int_{-1}^{\tau} T_k(\xi) d\xi  \\
= a_0\big(T_1(\tau) + T_0(\tau)\big)  + a_1\left( \frac{T_2(\tau)} {2(2)}- \frac{T_0 (\tau)}{2(2)} \right)\\
+ \sum_{k=2}^N a_k \left( \frac{T_{k+1}(\tau)}{2(k+1)} - \frac{T_{k-1}(\tau)}{2(k-1)}   \right)
- \sum_{k=2}^N a_k \left( \frac{(-1)^k T_0(\tau)}{k^2-1} \right)
\end{multline}
Rearranging the terms from the right-hand-side of the second equality in \eqref{eq:BVt-2}, we get,
\begin{multline}\label{eq:BVt-3}
\sum_{j=0}^N v_j B_j^a(\tau) =  T_0(\tau)\left( a_0  - \frac{a_1}{2(2)}  - \sum_{k=2}^N a_k \left( \frac{(-1)^k }{k^2-1} \right)    \right)\\
+ T_1(\tau)\left( a_0  - \frac{a_2}{2} \right) + \sum_{k=2}^{N-1} T_k(\tau)\left( \frac{a_{k-1} -a_{k+1}}{2k }  \right)  \\
+ T_N(\tau)\left( \frac{a_{N-1}}{2N} \right) + T_{N+1}(\tau)\left( \frac{a_{N}}{2(N+1)} \right)
\end{multline}
Thus, \eqref{eq:BVt-3} is equivalent to
\begin{equation}\label{eq:BVt-4}
\sum_{j=0}^N v_j B_j^a(\tau)  = \sum_{k=0}^{N+1} c_k T_k(\tau)
\end{equation}
where $c_k, k =0, \ldots, N+1$ are given by,
\begin{subequations}\label{eq:c_k=bydef}
\begin{align}
c_0 &:=  a_0  - \frac{a_1}{2(2)}  - \sum_{k=2}^N a_k \left( \frac{(-1)^k }{k^2-1} \right)
& c_1 &:= a_0  - \frac{a_2}{2}\\
c_k &:= \frac{a_{k-1} -a_{k+1}}{2k } & k & =2, \ldots, N-1\\
c_N &:= \frac{a_{N-1}}{2N}
& c_{N+1} &:= \frac{a_{N}}{2(N+1)}
\end{align}
\end{subequations}
Note that we now have an extra term with $c_{N+1}$.  That is, it is sufficient for \eqref{eq:BVt-4} to hold with $k$ running from $0$ to $N$:
\begin{equation}\label{eq:BVt-5}
\sum_{j=0}^N v_j B_j^a(\tau)  = \sum_{k=0}^{N} \widehat{c}_k T_k(\tau)
\end{equation}
The ``discrepancy'' between \eqref{eq:BVt-4} and \eqref{eq:BVt-5} is easily managed through the use of the Chebyshev aliasing formula; see Appendix C.  From the Chebyshev aliasing formula, the values of the $(N+1)^{th}$ Chebyshev polynomial, $T_{N+1}(\tau)$, evaluated at the $(N+1)$ CGL points, $\tau_0, \ldots, \tau_N$ are exactly the same as that of $T_{N-1}(\tau)$.  That is, $T_{N+1}(\tau)$ is an alias of $T_{N-1}(\tau)$ with $c_{N+1}T_{N+1}(\tau)$ and $c_{N-1}T_{N-1}(\tau)$ containing the same information as $\widehat{c}_{N-1}T_{N-1}(\tau)$. Hence, the coefficients in \eqref{eq:BVt-5} are given by,
\begin{subequations}\label{eq:chatk=}
\begin{align}
\widehat{c}_k &= c_k \quad \text{for } k = 0, \ldots, N-2, N \\
\widehat{c}_{N-1} &= c_{N-1} + c_{N+1}
\end{align}
\end{subequations}
It is straightforward to verify that computing $\widehat{c}_k, k =0, \ldots, N$ requires only $4N$ operations. Thus we have the following fast algorithm for computing $\bB^a V$:\\[1em]
%=========================================================
\fbox{%
    \parbox{0.95\linewidth}{%
\textbf{\texttt{FastBV}: A fast algorithm for computing the product $\bB^aV$}:
\begin{enumerate}
\item Given $v_0, v_1, \ldots, v_N$, apply \texttt{\textbf{Fastv2a}} (see Appendix B) to generate $a_0, a_1, \ldots, a_N$.
\item Compute $\widehat{c}_0, \widehat{c}_1, \ldots, \widehat{c}_N$ using \eqref{eq:chatk=}.
\item Apply \texttt{\textbf{Fasta2v}} to $\widehat{c}_0, \widehat{c}_1, \ldots, \widehat{c}_N$ to produce $\bB^aV$.
\end{enumerate}
}}\\ \\
%===============================================
%
It is obvious that the time complexity of \texttt{\textbf{FastBV}} is $\mathcal{O}(N\log(N)) + 4N + \mathcal{O}(N\log(N)) = \mathcal{O}(N\log(N))$ operations.  This point is demonstrated in Figure~\ref{fig:BVisNlogN}.
%
%======================================================================================
%==========================================================================================
%
The run time (ordinate) indicated in Figure~\ref{fig:BVisNlogN} is normalized to the computational time of $N= 2^8 = 256$.  The normalization is performed because, as noted in Reference~[\citen{newBirk-part-I}], actual run times are quite useless as a mathematical measure of computational efficiency.  Furthermore, to remove other practical factors in compute time\cite{newBirk-part-I}, the run times were averaged over 100 applications of the operations.

An important point to note in Figure~\ref{fig:BVisNlogN} is that for $N \lesssim 1,000$ ordinary matrix multiplication is faster than using the DCT. In other words, the advantage of using the fast DCT is more pronounced as $N$ increases past~$2^{10} = 1,024$. This point is theoretically consistent with computational complexity estimates that are valid asymptotically.

\section{Construction of a Low-Memory Fast Million-Point Linear Solver}
As implied in the Introduction section, it is very tempting to solve Problem~$(P^N_a)$ using an off-the-shelf optimizer like SNOPT or IPOPT. Such generic optimization routines do not allow matrix-free operations. Furthermore, generic optimizers are incapable of taking the added advantages of additional structural properties offered by Birkhoff-theoretic methods\cite{ross:Hessians}. In addition, it will be apparent shortly that to produce a fast Birkhoff-theoretic trajectory optimization solver, we need to construct a special preconditioner that needs to be inserted deep inside an optimization routine.  Hence, it is necessary to construct a new algorithm from the ground up.

In view of the preceding arguments, consider the problem of finding a feasible solution to Problem~$(P^N_a)$. To emphasize the Birkhoff elements of a new algorithm, we tentatively set $x^a$ and $x^b$ to some fixed value.   Then, it follows that we need an efficient method for solving the system of nonlinear equations given by,
\begin{subequations}\label{eq:X-V=nonlin}
\begin{align}
X  - \bB^a V  &= x^a\,\bb \\
V - f(X, U) & = 0
\end{align}
\end{subequations}
In anticipation of using a Newton's method, consider the derivatives (Jacobian) in \eqref{eq:X-V=nonlin}. The Jacobian of $f$ with respect to $X$ is given by,
\begin{multline}\label{eq:diag-partialf}
\partial_X f(X, U)  =
\left[
\begin{array}{ccc}
\partial_xf(x_0, u_0) & 0 & 0 \\
0 & \ddots & 0 \\
0 & 0 & \partial_xf(x_N, u_N) \\
\end{array}
\right] \\
 := diag[\partial_xf(x_0, u_0), \ldots, \partial_xf(x_N, u_N)]^T
\end{multline}
Because $V$ is linear in \eqref{eq:X-V=nonlin}, it follows that for fixed values of $u_0, \ldots, u_N$, we need an efficient method to solve the Newton system,
\begin{equation}\label{eq:Newton-ivp}
[I  - \bB^a \partial_X f]\delta X  = x^a\,\bb
\end{equation}
where $\delta X$ is the Newton step.
In \eqref{eq:Newton-ivp}, we have dropped the dependence of $X$ and $U$ in $\partial_X f$ because the Jacobian does not vary while solving \eqref{eq:Newton-ivp}.  Conceptually, a solution to \eqref{eq:Newton-ivp} is given by,
\begin{equation}\label{eq:Newton-ivp-inv}
\delta X  =  x^a\,[I  - \bB^a \partial_X f]^{-1}\bb
\end{equation}
Obviously, the worst algorithm to solve \eqref{eq:Newton-ivp} is to compute the matrix inverse first and multiply it by $\bb$ as implied by \eqref{eq:Newton-ivp-inv}.  A less naive approach to solving \eqref{eq:Newton-ivp} is to use a ``direct method'' such as an LU or QR decomposition method as is commonly used in many optimization routines\cite{snopt.paper,ipopt}. The computational complexity of these methods\cite{trefethen-bau-1997,golub-book-2013} is approximately $\mathcal{O}(N^3)$.  An obvious problem with this approach is that all the advantages of using an FFT to compute a Birkhoff matrix-vector product would be lost if a direct method were used to solve \eqref{eq:Newton-ivp}.  Hence a more efficient path to solving \eqref{eq:Newton-ivp} is to use a matrix-free method described earlier as Algorithm~1. Mapping the notation of \eqref{eq:Newton-ivp} to \eqref{eq:Ax=b} we have,
\begin{equation}\label{eq:Achibetadef}
A = [I  - \bB^a \partial_X f], \quad  \chi = \delta X  \quad \beta =  x^a\,\bb
\end{equation}
To use a matrix-free method for solving \eqref{eq:Newton-ivp} we first need to compute the product $A \chi$ without computing~$A$ first. This is achieved by the algorithm \textbf{\texttt{FastAX}}:\\ \\
%
%===============================================
\fbox{%
    \parbox{0.95\linewidth}{%
\textbf{\texttt{FastAX}: Matrix-free computation of $A\chi$ for \eqref{eq:Achibetadef} }:
\begin{enumerate}
\item Compute $[\partial_X f]\chi$.  The cost of this operation is only $\mathcal{O}(N)$ (because $[\partial_X f]$ is a diagonal matrix; see \eqref{eq:diag-partialf}).
\item Using the result of Step~1, compute $\bB^a \big([\partial_X f] \chi\big)$ using the algorithm \textbf{\texttt{FastBV}} discussed in the previous section.  The computational cost of this step is $\mathcal{O}(N\log(N))$.
\item Compute $\left(\chi - \bB^a \big([\partial_X f] \chi\big) \right)$.  The computational cost of this step is zero.
\end{enumerate}
}%
}
%}}
\\[1em]
%------------
Using \textbf{\texttt{FastAX}}, we can compute $A\chi$ in $\mathcal{O}(N\log(N))$ time (versus $\mathcal{O}(N^2)$ time if $\bB^a$ were computed first).  Hence, if an iterative method for solving \eqref{eq:Ax=b} is used, it is clear that we can compute $\delta X$ faster than $\mathcal{O}(N^3)$ provided the number of iterations are not too large. As is well-known\cite{CTKelley.1999,golub-book-2013,trefethen-bau-1997} if $A$ is not well-conditioned then then the number of these iterations become excessively large.  Although $[I-\bB^a]$ is well-conditioned (see Figure~\ref{fig:BirkCondNumbers}) it is obvious that the condition number of $A = [I  - \bB^a \partial_X f]$ depends strongly on the entries of $\partial_X f$.  In other words, the condition number of $A$ is dependent on the problem.  To overcome this problem-dependent computational issue, it is critical to design a preconditioner that works for arbitrary entries of $\partial_X f$ so that \emph{the goal of solving a million point trajectory optimization problem is not restricted to a small class of problems}.

\subsection{Fundamentals of Solving a Preconditioned System }
A generic preconditioner is based on the idea of producing a matrix $P$ whose inverse $P^{-1}$ is ``inexpensive to compute'' (e.g., when $P$ is a diagonal matrix) so that \eqref{eq:Ax=b} can be ``preconditioned'' by,
\begin{equation}\label{eq:leftP}
P^{-1}A\,\chi = P^{-1}\beta
\end{equation}
The second requirement for $P$ is that the product $(P^{-1}A)$ be better conditioned than $A$. Obviously, the best choice for $P$ is $A$ itself but it requires $A^{-1}$ to be inexpensively computable in which case \eqref{eq:leftP} becomes redundant.  Hence, a good choice for $P$ is one where $P \approx A$ with $P^{-1}$  inexpensively computable.  One issue in using \eqref{eq:leftP} is that the residual is also multiplied by $P^{-1}$ which distorts the accuracy of the computed solution as measured by $\norm{A\chi-\beta} \le \epsilon$.  An easy way to maintain the accuracy requirement is to \textit{right-precondition} \eqref{eq:Ax=b},
\begin{equation}\label{eq:rightP}
A\,P^{-1} (P \chi) = \beta
\end{equation}
Equation~(\ref{eq:leftP}) is called \textit{left preconditioning}. In right preconditioning, \eqref{eq:Ax=b} is solved by solving
\begin{subequations}\label{eq:rightPx2}
\begin{align}
A\,P^{-1} \Upsilon = \beta \label{eq:PinvUpsilon}\\
P \chi = \Upsilon \label{eq:Pchi}
\end{align}
\end{subequations}
where \eqref{eq:PinvUpsilon} is solved up to some tolerance $\epsilon > 0$.
Thus, Algorithm~1 for solving \eqref{eq:Ax=b} described in the Introduction section is modified to Algorithm~2:\\ \\
%
%===============================================
\fbox{%
    \parbox{0.95\linewidth}{%
\textbf{Algorithm 2: Rudimentary right-preconditioned matrix-free method for solving \eqref{eq:Ax=b}}:
\begin{enumerate}
\item Let $\chi^0$ be a guess to the solution of \eqref{eq:Ax=b}.
\item Compute the product $A\chi^0$ (without computing $A$ explicitly; see \textbf{\texttt{FastAX}}).
\item If $\norm{A\chi^0-\beta} \le \epsilon$, where $\epsilon > 0$ is a given solution tolerance, stop and exit.  Otherwise, compute $P\chi^0 = \Upsilon^0$ and continue.
\item For $k = 0, 1, \ldots$ compute $\Upsilon^{k+1} = \Upsilon^k + \Delta \Upsilon^k$ where $\Delta \Upsilon^k$ is found iteratively using matrix-vector products $AP^{-1}\Upsilon^k$ (without explicitly using $AP^{-1}$ ).
\item Stop the iteration when $\norm{AP^{-1}\Upsilon^k-\beta} \le \epsilon$.  Compute $\chi^k = P^{-1}\Upsilon^k$ and exit.
\end{enumerate}
}%
}
%}}
\\[1em]
%------------
\subsection{Development of a Special Birkhoff-Centric Preconditioner }
To find a good preconditioner for \eqref{eq:Newton-ivp}, we need to find a matrix $P$ that is close to
$[I  - \bB^a \partial_X f]$ and easily invertible. To motivate the development of such a matrix, we first examine $\bB^a$ and observe that it is lower left dominant as illustrated by the following matrix for $10$ CGL points:
%--------------------------
\begin{equation}\label{eq:Ba=10}
\bB^a_{10} = \left[
  \begin{array}{rrrrrrrrrr}
    0.00  &       0.00    &     0.00   &      0.00    &     0.00     &    0.00    &     0.00   &      0.00  &       0.00   &      0.00 \\
0.02  &  0.04  &  -0.00 &   0.00 &  -0.00 &   0.00 &  -0.00 &   0.00 &  -0.00  &   0.00 \\
    0.00  &  0.14  &  0.10  & -0.02  &  0.01 &  -0.01 &   0.00 &  -0.00  &  0.00   & -0.00 \\
    0.02  & 0.10  &  0.26   & 0.14 &  -0.02 &   0.01 &  -0.00 &   0.00   & -0.00    & 0.00 \\
    0.01   & 0.13  &  0.21  &  0.33  &  0.17 &  -0.03  &  0.02 &  -0.01 &   0.01   & -0.00 \\
    0.02  &  0.11  &  0.24   & 0.28 &   0.37  &  0.17 &  -0.03  &  0.02 &   -0.02    & 0.01 \\
    0.01  & 0.12   &  0.22  &  0.31  &  0.33 &   0.37 &   0.16 &  -0.03 &   0.02   & -0.01 \\
    0.01  &  0.11  &  0.23    &0.30 &   0.35   & 0.34  &  0.32  &  0.13&   -0.03  &  0.01 \\
    0.01  &  0.12  &  0.22   & 0.30  &  0.34  &  0.34  &  0.30 &   0.23 &   0.08 & -0.01 \\
    0.01  &  0.12  &  0.23  &  0.30 &   0.34 &   0.34 &   0.30  &  0.23  &  0.12 & 0.01 \\
\end{array}
\right]
\end{equation}
That is, the upper triangular elements of the Birkhoff-Chebyshev matrix are quite small relative to the lower triangular elements. Note that we have rounded the values of $\bB^a_{10}$ in \eqref{eq:Ba=10} to two decimal places.  In actual computation, $\bB^a$ must be computed to near machine precision. In any case, the last row of $\bB^a$ is exactly equal to the Clenshaw-Curtis weights by definition (see [\citen{newBirk-2023}]).  Now observe that each column of the lower triangular portion of $\bB^a_{10}$ is approximately a constant (excepting the diagonal elements). Furthermore, the diagonal elements are approximately half of the corresponding Clenshaw-Curtis weight. Hence, we can easily find an approximate $\bB^a$ matrix using the Clenshaw-Curtis weights $w^{CC}_j, \ j = 0, \ldots, N$ and the formula,
\begin{equation}
[\widetilde{\bB^a}]_{ij} := \left\{
                              \begin{array}{ll}
                                0, & \hbox{if } i < j \\
                                \displaystyle\frac{w^{CC}_j}{2} , & \hbox{if } i = j \\[0.5em]
                                w^{CC}_j, & \hbox{if } i > j
                              \end{array}
                            \right.
\end{equation}
The Clenshaw-Curtis weights can be computed at $\mathcal{O}(N\log(N))$ speed.  Obviously, $\widetilde{\bB^a}$ is lower triangular.  A $10\times10$ version of $\widetilde{\bB^a}$ is given by,
%--------------------------
\begin{equation}\label{eq:Ba=10}
\widetilde{\bB^a}_{10} = \left[
  \begin{array}{rrrrrrrrrr}
    0.01  &       0.00    &     0.00   &      0.00    &     0.00     &    0.00    &     0.00   &      0.00  &       0.00   &      0.00 \\
0.01  &  0.06  &  0.00 &   0.00 &  0.00 &   0.00 &  0.00 &   0.00 &  0.00  &   0.00 \\
    0.01  &  0.12  &  0.11  & 0.00  &  0.00 &  0.00 &   0.00 &  0.00  &  0.00   & 0.00 \\
    0.01  & 0.12  &  0.23   & 0.15 &  0.00 &   0.00 &  0.00 &   0.00   & 0.00    & 0.00 \\
    0.01   & 0.12  &  0.23  &  0.30  &  0.17 &  0.00  &  0.00 &  0.00 &   0.00   & 0.00 \\
    0.01  &  0.12  &  0.23   & 0.30 &   0.34  &  0.17 &  0.00  &  0.00 &   0.00    & 0.00 \\
    0.01  &  0.12   &  0.23  &  0.30  &  0.34 &   0.34 &   0.15 &  0.00 &   0.00   & 0.00 \\
    0.01  &  0.12  &  0.23    &0.30 &   0.34   & 0.34  &  0.30  &  0.11&   0.00  &  0.00 \\
    0.01  &  0.12  &  0.23   & 0.30  &  0.34  &  0.34  &  0.30 &   0.23 &   0.06 & 0.00 \\
    0.01  &  0.12  &  0.23  &  0.30 &   0.34 &   0.34 &   0.30  &  0.23  &  0.12 & 0.01 \\
\end{array}
\right]
\end{equation}
Note that $\widetilde{\bB^a}$ is more structured than simply being lower triangular. By observation, this structure can be used advantageously to compute $\widetilde{\bB^a} V$ using the recursion,
\begin{equation}
R_k = R_{k-1} - \frac{w_{k-1}v_{k-1}}{2} + w_{k-1}v_{k-1} + \frac{w_k v_k}{2}, \quad  k = 1, \ldots, N, \quad R_0 = \frac{w_0 v_0}{2}
\end{equation}
where $R_k$ is the $k^{th}$ row of $\widetilde{\bB^a}V = [\widetilde{\bB^a}V]_k $.  A more efficient $\mathcal{O}(N)$  algorithm for computing $\widetilde{\bB^a}V$ is described by the following algorithm:\\[1em]
%=========================================================
\fbox{%
    \parbox{0.95\linewidth}{%
\textbf{\texttt{FastBccV}: A fast algorithm for computing the product $\widetilde{\bB^a}V$}:
\begin{enumerate}
\item Compute $w^{CC}_j, \ j = 0, \ldots, N$ using an $\mathcal{O}(N\log(N))$ algorithm.
\item Set $s = 0$.
\item
   \begin{algorithmic}
     \For{$k = 0$ to $N$}
     \State $\xi_k = w^{CC}_k v_k/2 + s$
     \State $s \gets 2\xi_k - s$ \vskip 2mm
    \State $[\widetilde{\bB^a}V]_k = \xi_k$ \quad \% This line is written for clarity
     \EndFor
   \end{algorithmic}
\end{enumerate}
}}\\ \\ \\
%===============================================
%
Note that we never store or compute $\widetilde{\bB^a}$ explicitly to produce $\widetilde{\bB^a} V$. The computational space complexity for $\widetilde{\bB^a} V$ is just $\mathcal{O}(N)$.

All of these considerations lead us to define a special Birkhoff preconditioner according to,
\begin{equation}
P := [I  - \widetilde{\bB^a} \partial_X f]
\end{equation}

\subsection{A Preconditioned Algorithm for Solving \eqref{eq:Newton-ivp}}
\begin{figure}[h!]
  \centering
  \includegraphics[width=0.5\textwidth]{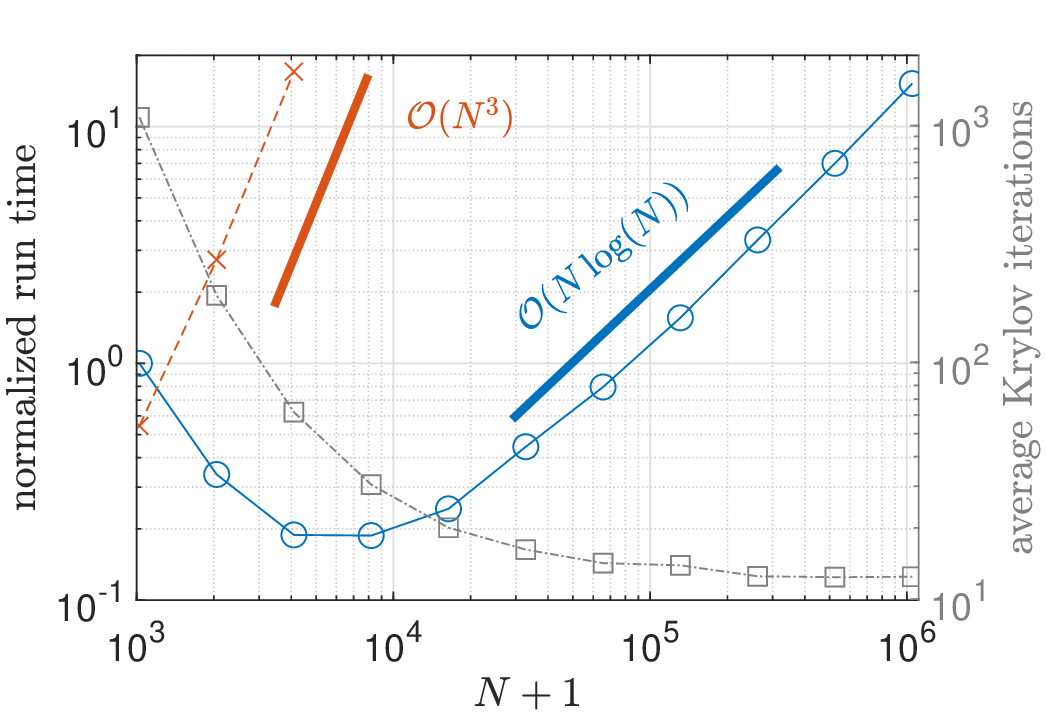}
  \caption{Normalized Run Times for \texttt{\textbf{FastLinSol}} for $N$ up to One Million Mesh/Grid Points. The $\mathcal{O}(N^3)$ Plot Uses MATLAB's Backslash Operation. }
  \label{fig:LinSol}
\end{figure}

Solving \eqref{eq:Newton-ivp} via \eqref{eq:rightPx2} requires a matrix-free computation of $A\,P^{-1} \Upsilon $ where $A$ is given by \eqref{eq:Achibetadef}.  To design an efficient algorithm to compute $P^{-1} \Upsilon$ observe that $\widetilde{\bB^a} \partial_X f$ has the same structure as $\widetilde{\bB^a}$ owing to the fact that $\partial_X f$ is diagonal. Furthermore, because $I$ is diagonal, the lower triangular portion of $P$ retains the same structure as that of $\widetilde{\bB^a}$. Hence $P^{-1} \Upsilon = [I  - \widetilde{\bB^a} \partial_X f]^{-1} \Upsilon $  can be computed by modifying the idea incorporated in \textbf{\texttt{FastBccV}}.  This algorithm is given by the following:\\[1em]
%=========================================================
\fbox{%
    \parbox{0.95\linewidth}{%
\textbf{\texttt{FastPinvY}: A fast algorithm for computing the product $P^{-1} \Upsilon = [I  - \widetilde{\bB^a} \partial_X f]^{-1} \Upsilon $ }
\begin{enumerate}
\item Compute $w^{CC}_k, \ k = 0, \ldots, N$ using an $\mathcal{O}(N\log(N))$ algorithm.
\item Compute $ \partial_{x_k}f, \ k = 0, \ldots, N$.  \quad \% This computation is part of the oracle complexity.
\item Set $s = 0$.
\item
   \begin{algorithmic}
     \For{$k = 0$ to $N$}
     \State $\xi_k = (1 - w^{CC}_k (\partial_{x_k}f)/2)^{-1}(\Upsilon_k - s)$
     \State $s \gets w^{CC}_k (\partial_{x_k}f) \Upsilon_k + s$ \vskip 2mm
    \State $[P^{-1} \Upsilon]_k = \xi_k$ \quad \% This line is written for clarity
     \EndFor
   \end{algorithmic}
\end{enumerate}
}}\\ \\
%===============================================
%
Combining \textbf{\texttt{FastPinvY}} with \textbf{\texttt{FastAX}} primes us to compute $[AP^{-1}\Upsilon]$ at  $\mathcal{O}(N\log(N))$ compute time for insertion into Algorithm~2 to solve \eqref{eq:Newton-ivp}.  Collecting all of the preceding ideas, the details of this algorithm are as follows:\\ \\
%
%===============================================
\fbox{%
    \parbox{0.95\linewidth}{%
\textbf{\texttt{FastLinSol}: Right-preconditioned matrix-free method for solving \eqref{eq:Newton-ivp}}:
\begin{enumerate}
\item Let $\delta X^0$ be a guess to the solution of \eqref{eq:Newton-ivp}.
\item Compute the product $[I  - \bB^a \partial_X f]\delta X^0$ using the fast matrix-free algorithm \textbf{\texttt{FastAX}}.
\item If $\norm{[I  - \bB^a \partial_X f]\delta X^0 -x^a\,\bb} \le \epsilon$, where $\epsilon > 0$ is a given solution tolerance, stop and exit.  Otherwise, compute $\Upsilon^0 \leftarrow [I  - \widetilde{\bB^a} \partial_X f]\delta X^0$ using \textbf{\texttt{FastAX}} (modified by \texttt{\textbf{FastBccV}}) and continue.
\item For $k = 0, 1, \ldots$ compute $\Upsilon^{k+1} = \Upsilon^k + \Delta \Upsilon^k$ where $\Delta \Upsilon^k$ is found iteratively using \textbf{\texttt{FastPinvY}} and \textbf{\texttt{FastAX}}. See \eqref{eq:rightPx2}.
\item Stop the iteration when $\norm{[I  - \bB^a \partial_X f][I  - \widetilde{\bB^a} \partial_X f]^{-1}\Upsilon^k-x^a\,\bb} \le \epsilon$.  Compute $\delta X^k$ using \textbf{\texttt{FastPinvY}} and exit.
\end{enumerate}
}%
}
%}}
\\[1em]
%------------
%
Figure~\ref{fig:LinSol} shows the results of a scaling experiment using \texttt{\textbf{FastLinSol}} using a GMRES implementation for Step~4.
%
%===============================================
%===================================
%
The run times in Figure~\ref{fig:LinSol} are normalized to run times for $N = $1,024.  An interesting feature of Figure~\ref{fig:LinSol} is that the runtime initially produces \emph{faster} solves for \emph{larger} values of $N$.  This is likely because for smaller values of $N$, $\widetilde{\bB^a}$ is further away from $\bB^a$ and that as $N$ increases $\widetilde{\bB^a}$ approaches $\bB^a$ making  $AP^{-1}$ in \eqref{eq:PinvUpsilon} closer to the identity matrix.  This hypothesis is verified by observing the number of Krylov iterations shown in Figure~\ref{fig:LinSol}.
Also shown in Figure~\ref{fig:LinSol} is the (normalized) run time when \eqref{eq:Newton-ivp} is solved using MATLAB's backslash operator. In this setting the backslash operator is a call to UMFPACK \cite{mldivide.2024,davis.2004}. It is apparent that the backslash operator is substantially slower and scales as $\mathcal{O}(N^3)$.

\section{An Illustrative Nonlinear Nonconvex Astrodynamics Problem}
\begin{figure}[h!]
  \centering
  \includegraphics[width=0.45\textwidth]{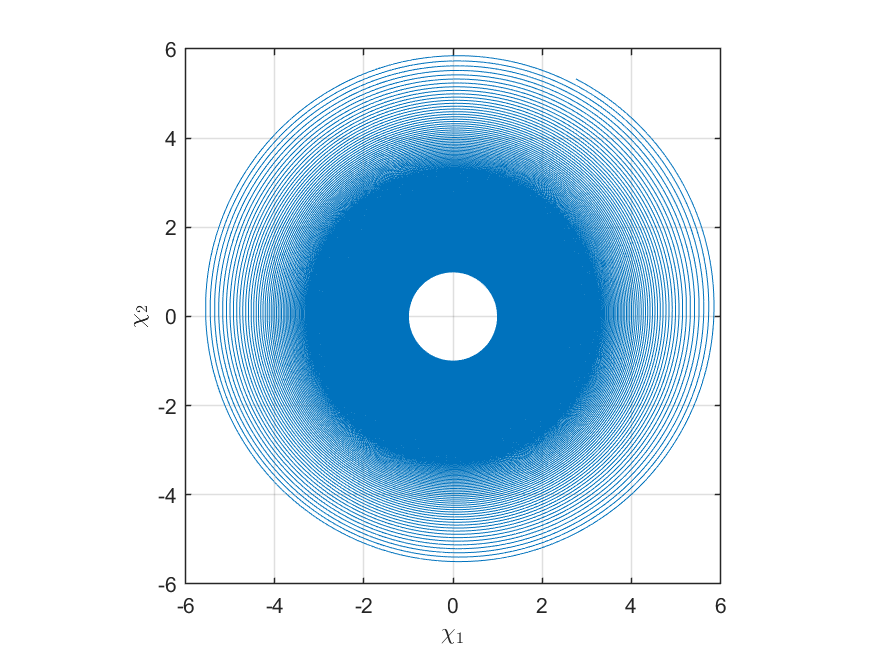}
  \includegraphics[width=0.45\textwidth]{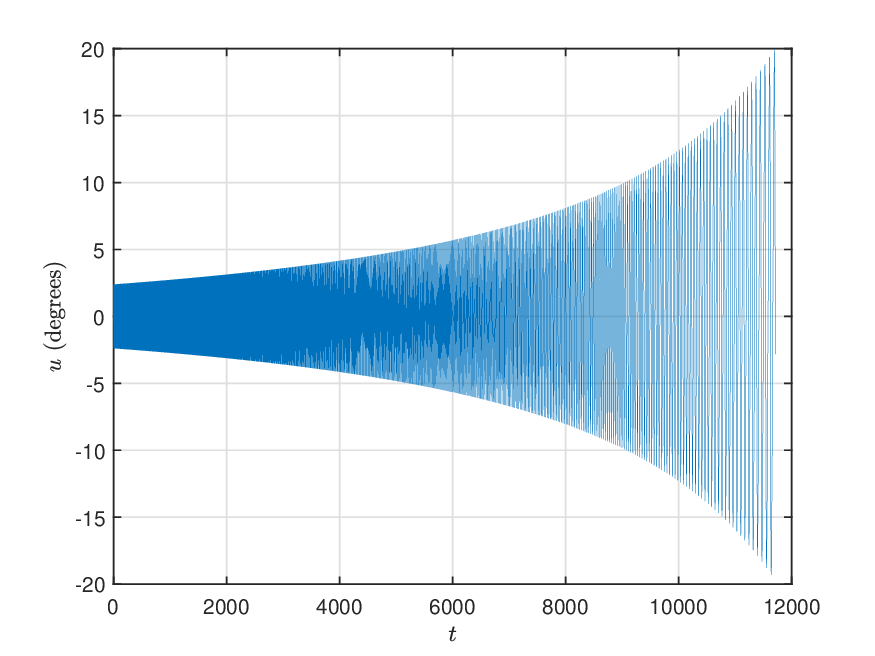}
  \caption{Candidate optimal trajectory and control solution for the orbit transfer problem (\ref{eq:ot}) with $N+1=$1,048,577 CGL points. }
  \label{fig:ot}
\end{figure}

Consider the low-thrust orbit transfer problem discussed in Reference~[\citen{newBirk-part-II}].
This problem is given by,
\begin{subequations}\label{eq:ot}
\begin{align}
  \underset{y(\cdot),u(\cdot),T}{\text{minimize}}&\quad J[y(\cdot),u(\cdot),T]:= T \\
    \text{subject to}&\quad \frac{dy}{dt} = \displaystyle\begin{cases}
  y_1\\
  \frac{y_2^2}{y_0} - \frac{1}{y_0^2} + a \sin u\\
  -\frac{y_1y_2}{y_0} + a \cos u\\
  \frac{y_2}{y_0}
  \end{cases}\\
  \text{and}&\quad y(0) = (1,0,1, 0)^\top, \quad\text{ }\quad y(T) = (6,0,1/\sqrt{6})^\top
\end{align}
\end{subequations}
%======================================
A Birkhoff-theoretic solution to this problem generated over a million CGL points is shown in Figure~\ref{fig:ot}.
%
%--------------------------------------------------------------
%-----------------------------------------------------
%
This solution was verified and validated using Pontryagin's Principle\cite{ross-book} using the method described in Reference~[\citen{newBirk-part-II}].  The first point to note in  Figure~\ref{fig:ot} is the production of a solution to a practical, nonlinear, nonconvex astrodynamics problem over a million grid points $(= 2^{20}+1)$. The second point of note is that this solution was obtained on an ordinary laptop (MacBook Air) with a maximum RAM capacity of only 32GB. The computational speed for generating this solution was already shown in Figure~\ref{fig:LinSol} in terms of the major bottleneck in solving \eqref{eq:Newton-ivp}.  With regards to the total number of iterations, the statistics are shown in Table~\ref{tab:m}. From this table, it follows that the number of SQP iterations as well as the number solves of the KKT system are consistently the same across the 21 different values of $N$ that span three orders of magnitude.

%===============================================
\begingroup
\begin{table}[hbt]
  \setlength{\tabcolsep}{3pt}
  \centering
  \caption{Illustrating the mesh independence of solving a nonlinear, nonconvex astrodynamics problem.}
  \label{tab:m}
  \begin{tabular}{lllllllllllll}
\hline
    $N$                     && $2^{10}$ & $2^{11}$ & $2^{12}$ & $2^{13}$ & $2^{14}$ & $2^{15}$ & $2^{16}$ & $2^{17}$ & $2^{18}$ & $2^{19}$ & $2^{20}$ \\
    \hline
    SQP iterations          && 4 & 4 & 4 & 4 & 4 & 4 & 4 & 4 & 4 & 4 & 4 \\
    Matrix-free KKT solves  && 37 & 37 & 37 & 37 & 37 & 37 & 37 & 37 & 37 & 37 & 37\\
\hline
 \end{tabular}
\end{table}
\endgroup
%==============================================

\section{Conclusion}
To the best of the authors' knowledge, this is the first time a nontrivial, nonlinear, nonconvex trajectory optimization problem has been solved over a million grid points. Beyond this accomplishment, we have shown that such problems can be solved quickly and on ordinary present-day (2025) laptops.  The key to achieving this high performance on low-end processors is to use matrix-free methods. Matrix-free methods work most efficiently when the condition numbers are low.  This is why a Bikrhoff-theoretic method became a starting point to take advantage of the efficiency of matrix-free methods. The implementation of these ideas are far from complete from the point of view of generating production software.  However, we envision such fast solvers being available in the near future.

\section{Acknowledgment}
The work of the first three authors (AJ, DPK \& DR) was co-sponsored by the Laboratory Directed Research and Development program at Sandia National Laboratories.\footnote{
Sandia National Laboratories is a multi-mission laboratory managed and operated by National Technology \& Engineering Solutions of Sandia, LLC (NTESS), a wholly owned subsidiary of Honeywell International Inc., for the U.S. Department of Energy's National Nuclear Security Administration (DOE/NNSA) under contract DE-NA0003525. This written work is authored by an employee of NTESS. The employee, not NTESS, owns the right, title and interest in and to the written work and is responsible for its contents. Any subjective views or opinions that might be expressed in the written work do not necessarily represent the views of the U.S. Government. The publisher acknowledges that the U.S. Government retains a non-exclusive, paid-up, irrevocable, world-wide license to publish or reproduce the published form of this written work or allow others to do so, for U.S. Government purposes. The DOE will provide public access to results of federally sponsored research in accordance with the DOE Public Access Plan.} JDS was supported in part by NSF Grant DMS-2231482.  We also thank G. von Winckel for suggesting the aliasing formula in the \texttt{\textbf{FastBV}} algorithm.

\appendix

\section{APPENDIX A: CONNECTION BETWEEN CHEBYSHEV AND FOURIER SERIES}
The connection between Chebyshev and Fourier series is well-known\cite{boyd,trefethen-2000}. In fact, most books on spectral methods introduce a Fourier series as a precursor to a Chebyshev pseudospectral method.  The following is adapted from [\citenum{trefethen-2000}].

Let\footnote{To maintain consistency with standard notation, we use the symbol $i$ in this appendix to mean $\sqrt{-1}$.  Elsewhere, $i$ is an index, $0, 1, \ldots$.} $z = e^{i\theta}, \ i = \sqrt{-1}$ so that $\Re(z) = (z + z^{-1})/2$. This implies,
\begin{equation}\label{eq:zk=}
\Re(z^k) = (z^k + z^{-k})/2
\end{equation}
Next, consider rewriting $\Re(z^{k+2})$ as follows:
\begin{multline}
\Re\left(z^{(k+2)}\right) = \frac{1}{2} \left(z^{k+2} + z^{-(k+2)}\right)
 = \frac{1}{2} \left(z^{k+1}z + z^{-(k+1)}z^{-1}\right) \\
 =\left(\frac{z^{(k+1)} + z^{-(k+1)}}{2}\right)(z + z^{-1}) - \left(\frac{z^{(k)} + z^{-(k)}}{2}\right)
 = 2 \Re\left(z^{(k+1)}\right) \Re\left(z\right) - \Re\left(z^{(k)}\right) \label{eq:recursion-z}
%& = 2 \cos((k+1)\theta)\cos(\theta) - \cos(k\theta)
\end{multline}
Equation~(\ref{eq:recursion-z}) is identical to the Chebyshev recurrence relation,
\begin{equation}\label{eq:recurrence-T}
T_{k+2}(\tau) = 2 \tau T_{k+1}(\tau) - T_{k}(\tau), \quad k = 0, \ldots
\end{equation}
Hence, \eqref{eq:recursion-z} is equivalent \eqref{eq:recurrence-T} under the transformation,
\begin{multline}\label{eq:T=cos}
\tau = (z + z^{-1})/2 = \cos(\theta), \quad T_k(\tau) = T_k(\cos(\theta)) =  \Re(z^k) = (z^k + z^{-k})/2 \\
= \cos(k\theta(\tau)), \quad \theta(\tau) = \cos^{-1}(\tau)
\end{multline}
Note that \eqref{eq:T=cos} puts $\theta = 0$ at $\tau = 1$ and $\theta = \pi$ at $\tau = -1$.  See \eqref{eq:CGLpts} for context.

Consider next a Chebyshev series expansion of the virtual variable $v(\tau)$ (see \eqref{eq:xdot=v=sumT}),
\begin{equation}\label{eq:v=sumT}
v(\tau)  = \sum_{k=0}^N a_k T_k(\tau)  %= \sum_{j=0}^N v_j \dot B_j^a(\tau)
\end{equation}
Substituting \eqref{eq:T=cos} in \eqref{eq:v=sumT} we get,
\begin{equation}
v(\tau)  = \sum_{k=0}^N a_k \cos(k\theta(\tau))
\end{equation}
which is simply a Fourier series without the sine terms.

\section{APPENDIX B: AN FFT ALGORITHM FOR COMPUTING THE CHEBYSHEV MODAL COEFFICIENTS}
In this section, we lean heavily on [\citenum{newBirk-2023}], which in turn, uses many other sources.

Taking orthogonal projections of \eqref{eq:xdot=v=sumT} along the Chebyshev polynomials generates,
\begin{equation}\label{eq:ak=bydef}
a_k := \frac{1}{\gamma_k}\int_{-1}^1 \frac{v(\tau) T_k(\tau)}{\sqrt{1- \tau^2}} d\tau
\qquad \gamma_k := \int_{-1}^1 \frac{T_k(\tau) T_k(\tau)}{\sqrt{1- \tau^2}} d\tau
\end{equation}
where $\gamma_k$ is the normalization coefficient.
The normalization coefficient is given in closed form by
\begin{equation}
\gamma_k = \left\{
                   \begin{array}{ll}
                   \pi, & \hbox{if } k = 0 \\
                   \pi/2, & \hbox{if } k \ge 1
                   \end{array}
            \right.
\end{equation}
The integral pertaining to $a_k$ in \eqref{eq:ak=bydef} can also be evaluated in closed form using \eqref{eq:xdot=v=sumT}.  This formula is given by,
\begin{equation}
\int_{-1}^1 \frac{v(\tau) T_k(\tau)}{\sqrt{1- \tau^2}} d\tau = \left\{
                                                                 \begin{array}{ll}
                                                                    \displaystyle\sum_{j=0}^{N} v_j T_k(\tau_j) w_j^{CGL}, & \hbox{for  } k = 0, \ldots, N-1 \\
                                                                 \displaystyle\sum_{j=0}^{N} (v_j/2) T_k(\tau_j) w_j^{CGL}, & \hbox{for } k = N
                                                                 \end{array}
                                                               \right.
\end{equation}
where $w_j^{CGL}$ are the CGL quadrature weights given by $w_0^{CGL} = \pi/(2N) = w_N^{CGL}, w_j^{CGL} = \pi/N, j = 1, \ldots, N-1$.  Combining all these results, it is straightforward to show that $a_k, k = 0, \ldots, N$ simplify to,
\begin{subequations}\label{eq:ak=3}
\begin{align}
%a_0 &= \frac{1}{N} \left(\frac{v_0}{2} + \frac{v_N}{2} + \sum_{j=1}^{N-1}v_j  \right)\\
a_k & = \frac{2}{N} \left(\frac{v_0}{2}(-1)^k + \frac{v_N}{2} + \sum_{j=1}^{N-1}v_j  \cos\left(\frac{k(N-j) \pi}{N} \right) \right), \quad k = 1, \ldots, N-1\\
a_0 &= \frac{1}{N} \left(\frac{v_0}{2} + \frac{v_N}{2} + \sum_{j=1}^{N-1}v_j  \right) \quad a_N  = \frac{1}{N} \left(\frac{v_0}{2} (-1)^N + \frac{v_N}{2} + \sum_{j=1}^{N-1}v_j \cos\big((N-j)\pi\big) \right)
\end{align}
\end{subequations}
To connect it to DCT-I (see \eqref{eq:vm=DCT-1}), we define
\begin{equation}
\widetilde{a}_0 = 2 a_0, \quad \widetilde{a}_N = 2 a_N, \quad \widetilde{a}_k = a_k
\end{equation}
Then, \eqref{eq:ak=3} can be rewritten as,
\begin{align}
\widetilde{a}_k &= \frac{2}{N} \left( \frac{v_0}{2}(-1)^k  +   \frac{v_N}{2}
+ \sum_{j=1}^{N-1}v_j  \cos\left(\frac{k(N-j) \pi}{N} \right) \right),  \quad k = 0, \ldots, N
\end{align}
Thus, an FFT algorithm for computing the Chebyshev modal coefficients can be described as follows:\\[1em]
\fbox{%
    \parbox{0.95\linewidth}{%
\textbf{\texttt{Fastv2a}: A fast algorithm for computing the Chebyshev modal coefficients from its nodal values}:
\begin{enumerate}
\item Given $v_0, v_1, \ldots, v_N$, rearrange it to $v_N, v_{N-1}, \ldots, v_0$.
\item Apply DCT-I to $v_N, v_{N-1}, \ldots, v_0$ and multiply the result by $2/N$ to generate $\widetilde{a}_0, \widetilde{a}_1, \ldots, \widetilde{a}_N$.

\item Divide the first and last term of $\widetilde{a}_0, \widetilde{a}_1, \ldots, \widetilde{a}_N$ by two to produce $a_0, a_1, \ldots, a_N$.
\end{enumerate}
}}

\section{APPENDIX C: THE CHEBYSHEV ALIASING FORMULA}
The Chebyshev aliasing formula is well-known.\cite{boyd,trefethen-2000}  From \eqref{eq:T=cos} we have,
\begin{equation}\label{eq:boyd-1}
T_{N+1}(\cos\theta) = \cos[(N+1)\theta] = \cos(N\theta)\cos(\theta) - \sin(N\theta)\sin(\theta)
\end{equation}
Over the CGL grid (see \eqref{eq:CGLpts}) $\theta_j = j\pi/N, j = 0, \dots, N$.  As a result, $\sin(N\theta_j) = 0, j = 0, \ldots, N$. Substituting this result in \eqref{eq:boyd-1} we get,
\begin{equation}\label{eq:boyd-2}
T_{N+1}(\cos\theta_j) = \cos(N\theta_j)\cos(\theta_j), \quad j = 0, \ldots, N
\end{equation}
Similarly, we get,
\begin{equation}\label{eq:boyd-3}
T_{N-1}(\cos\theta_j) = \cos(N\theta_j)\cos(-\theta_j) = \cos(N\theta_j)\cos(\theta_j) =  T_{N+1}(\cos\theta_j), \quad j = 0, \ldots, N
\end{equation}
%

%From \eqref{eq:Tkj=coskj} it follows that,
%%
%\begin{equation}
%T_{N}(\tau_j) =  \cos\left((N-j) \pi\right), \quad j = 0, \ldots, N
%\end{equation}
%%
%
%$$T_{2lN \pm m}(\tau_j) =  \cos\left(\frac{(2lN \pm m ) (N-j) \pi}{N} \right), \quad j = 0, \ldots, N$$

%\bibliographystyle{AAS_publication}   % Number the references.
%\bibliography{references}   % Use references.bib to resolve the labels.

\end{document}